\documentclass[11pt]{amsart}

\usepackage{latexsym}
\usepackage{amsfonts}
\usepackage{amsthm}
\usepackage{graphicx}
\usepackage{amssymb}
\usepackage{amsmath}
\usepackage{amssymb}
\usepackage{color}
\author[]{Vladimir Shpilrain}
\address{Department of Mathematics, The City  College  of New York, New York, 
NY 10031}
\email{shpil@groups.sci.ccny.cuny.edu}

\author[]{Alexander Ushakov}
\address{Department of Mathematics, CUNY Graduate Center, New York, 
NY 10016}
\email{aushakov@mail.ru}
\thanks{Research of the first author was partially supported by 
the NSF grant DMS-0405105. Research of the 
second author was partially supported by Umbanet Inc. through an 
award from the U.S. Department of Commerce NIST, Advanced Technology 
Program, Cooperative Agreement No. 70NANB2H3012}

\begin{document}

\title[]{The conjugacy search problem in 
public key cryptography: unnecessary and insufficient}

\begin{abstract} 
The conjugacy search problem in a group $G$ is the problem 
of recovering an $x \in G$ from given $g \in G$ and $h=x^{-1}gx$. 
This problem is in the core of several recently suggested 
public key exchange protocols, most notably the one due to 
Anshel, Anshel, and Goldfeld, and the one due to Ko, Lee at al. 

 In this note, we make two observations that seem to have 
eluded most people's attention. The first observation 
is that solving the conjugacy search problem is not necessary 
for an adversary to get the common secret key in the Ko-Lee 
protocol. It is sufficient to solve an apparently easier problem 
of finding $x, y \in G$ such that $h=ygx$ for given $g, h \in G$. 

 Another observation 
is that solving the conjugacy search problem is not sufficient 
for an adversary to get the common secret key in the 
Anshel-Anshel-Goldfeld protocol. 
\end{abstract}

\maketitle

\section{Introduction}

One of the possible generalizations of the {\it discrete logarithm
problem} to arbitrary groups is the so-called
 {\it  conjugacy search  problem} (CSP): given two elements $g, h$ of a group $G$
 and the information that $g^x=h$
for some $x \in G$, find at least one particular element $x$ like
that. Here $g^x$ stands for $x^{-1}gx$. The (alleged)
computational difficulty of this  problem in some particular
groups (namely, in braid groups) has been  used in several group
based public key protocols, most notably in \cite{AAG}  and
\cite{KLCHKP}.

  In this note, we show that solving the conjugacy search  problem 
is unnecessary for an adversary to get the common secret key in the Ko-Lee 
(or any similar) protocol, and, on the other hand, 
 is insufficient to get the common secret key in the more sophisticated 
Anshel-Anshel-Goldfeld protocol. This raises the stock of the 
latter protocol and makes one think there might be more to it than 
meets the eye.

\section{Why solving CSP is unnecessary}

 First we   recall the (generalized) Ko-Lee protocol. A group $G$ 
(with efficiently  solvable word problem) 
and two commuting subsets $A, B \subseteq G$ (i.e., $ab=ba$ for any  $a \in A, 
 ~b \in B$) are public. An element  $w \in G$ is public, too. 

\begin{enumerate}
 \item[{\bf (1)}]  Alice selects a private $a \in A$ and sends the  
element $a^{-1} w a$ to Bob.

 \item[{\bf (2)}]  Bob selects a private $b \in B$ and sends the  
element $b^{-1} w b$ to Alice. 

 \item[{\bf (3)}] Alice computes $K_A=a^{-1} b^{-1} w b a$,
and Bob computes $K_B=b^{-1} a^{-1} w a b$. Since $a b=b a $ 
(and therefore, $a^{-1} b^{-1}=b^{-1}a^{-1}$) 
in $G$, one has  $K_A=K_B=K$ (as an element of $G$), which is now Alice's
 and Bob's common secret key.

\end{enumerate}

 Note that since we want the key space to be as big as possible, 
we may assume, to simplify the language in what follows, that, say,  
the set $A$ is maximal with the property that $ab=ba$ for any  $a \in A, 
 ~b \in B$.

 Now suppose an adversary finds $a_1, a_2$ such that 
$a_1 w a_2 = a^{-1} w a$ and $b_1, b_2$ such that 
$b_1 w b_2 = b^{-1} w b$. Suppose also that both $a_1, a_2$ commute with any 
$b \in B$. 
Then the adversary  gets 
$$a_1 b_1 w b_2 a_2 =
a_1 b^{-1} w b a_2 = b^{-1} a_1 w a_2 b = b^{-1} a^{-1} w a b =K.$$ 

 We emphasize that these $a_1, a_2$  and $b_1, b_2$ do not have to 
do anything with the private elements originally selected by Alice
 or Bob, which simplifies the search substantially. 

 In other words, to get the secret key $K$, the adversary does not have to solve
the conjugacy search  problem, but instead, it is sufficient to solve an 
apparently easier problem which some authors (see e.g. \cite{CKLHC})
call the {\it decomposition problem}: 

\medskip

 {\it Given an element $w$ of a group $G$ and another element
$x\cdot w\cdot y$, find any  elements $x'$ and $y'$ that would 
belong to a given subset $A \subseteq G$  and       satisfy  
$x'\cdot w\cdot y' = x\cdot w\cdot y$.}

\medskip

 We note that the   condition $x', y' \in A$ 
may not be easy to verify for some subsets $A$, but for the 
particular situation considered in \cite{KLCHKP} this is 
straightforward and can be done just by inspection of the 
normal forms of $x$ and  $y$. 

 The claim that the decomposition problem should be easier 
than the conjugacy search  problem is intuitively 
clear since it is generally  easier to solve an equation 
with two unknowns than a special case of the same equation 
with just one unknown.

\section{Why solving CSP is insufficient}

 The protocol that we describe below, due to Anshel, Anshel, and Goldfeld 
\cite{AAG},  
is more complex than the protocol in the previous section, 
but it is more general in the sense that there are no requirements  
on the group $G$ other 
than to have efficiently  solvable word problem. This really makes a difference 
and gives a big advantage to the protocol of \cite{AAG} over that 
of \cite{KLCHKP}.

 A group $G$ and   elements $a_1,...,a_k, b_1,...,b_m \in G$ are public. 

\begin{enumerate}
 \item[{\bf (1)}] Alice picks a private $x \in G$ as a word in 
$a_1,...,a_k$ (i.e.,  $x=x(a_1,...,a_k)$)  and sends 
$b_1^x,...,b_m^x$ to Bob. 

 \item[{\bf (2)}] Bob picks a private $y \in G$ 
as a word in $b_1,...,b_m$ and sends 
$a_1^y,...,a_k^y$ to Alice.

 \item[{\bf (3)}] Alice computes $x(a_1^y,...,a_k^y) = x^y = y^{-1}xy$, 
and Bob computes \\ 
$y(b_1^x,...,b_m^x) = y^x = x^{-1}yx$. Alice and Bob then come up with a 
common private key $K=x^{-1}y^{-1}xy$ 
(called the \emph{commutator} of $x$ and $y$) as follows:  Alice multiplies 
$y^{-1}xy$ by $x^{-1}$ on the left,  while
Bob multiplies 
$x^{-1}yx$ by $y^{-1}$ on the left, and then takes the inverse of 
the whole thing: $(y^{-1}x^{-1}yx)^{-1} = x^{-1}y^{-1}xy$.

\end{enumerate}

 It appears to be a common belief (see e.g. \cite{GKTTV1,  HS, HT}) 
 that solving the conjugacy search problem 
for $b_1^x,...,b_m^x, a_1^y,...,a_k^y$ in the group $G$ would allow   an adversary  
to get the secret key $K$. However, if we look at Step (3) of the 
protocol, we see that the adversary would have to know, say,  $x$ not 
simply as  a word in the generators of the group $G$, but as  a word in 
$a_1,...,a_k$. That means the adversary would also have to solve 
the {\it membership search problem}: 

\medskip

 {\it Given   elements $x, a_1,...,a_k$ of a group $G$, 
find an expression (if it exists) of  $x$ as a word in 
$a_1,...,a_k$.} 
\medskip

 We note  that the (decision)  membership problem is to 
determine whether or not a given $x \in G$ belongs to the subgroup of 
$G$ generated by given $a_1,...,a_k$. Even this, apparently 
easier problem, turns out to be quite hard in most groups. 
For instance, the membership problem in a braid group $B_n$ is
algorithmically unsolvable if $n \ge 6$ because such a braid group contains
 subgroups isomorphic to  $F_2 \times F_2$ (that would 
be, for example, the subgroup generated by $\sigma_1^2,
\sigma_2^2, \sigma_4^2$, and  $\sigma_5^2$, see \cite{Collins}), where $F_2$ 
is the free group of rank 2. In the group $F_2 \times F_2$, 
the membership problem is algorithmically unsolvable by 
an old result of Mihailova \cite{Mihailova}. 

 We also note that if the adversary finds, say, some $x' \in G$ 
such that ~$b_1^x=b_1^{x'}, ... , \\
b_m^x=b_m^{x'}$, there is no guarantee 
that $x'=x$ in $G$.  Indeed, if $x'=c_b x$, where $c_b b_i= b_i c_b$ for all $i$, 
then $b_i^x=b_i^{x'}$ for all $i$, and therefore  $b^x=b^{x'}$ for any 
element $b$  from the subgroup  generated by $b_1,...,b_m$; in particular, 
$y^x=y^{x'}$. Now the problem is that if $x'$  does not belong to the 
subgroup $A$ generated by $a_1,...,a_k$ (which may very well be the case), 
 then the adversary   will not be able to obtain the common secret key $K$.
On the other hand, if  $x'$ (and, similarly, $y'$)  does   belong to the 
subgroup $A$ (respectively, to the subgroup $B$ generated by $b_1,...,b_m$), then 
the adversary   will  be able to get the correct $K$ even though his 
$x'$ and $y'$ may be different from $x$ and $y$, respectively. 
Indeed, if $x'=c_b x$, $y'=c_a y$, where $c_b$ centralizes $B$  and 
 $c_a$ centralizes $A$, then 
$$x'^{-1}y'^{-1}x'y'=(c_b x)^{-1}(c_a y)^{-1}c_b x c_a y=x^{-1}c_b^{-1}y^{-1}c_a^{-1}c_b x c_a y=
x^{-1}y^{-1}xy=K$$
 because $c_b$ commutes with $y$ and with  $c_a$ (note that $c_a$ belongs 
to the subgroup $B$, which follows from the assumption $y'=c_a y \in B$, 
and, similarly, $c_b$ belongs to $A$), and $c_a$ commutes with $x$. 

We emphasize that the adversary ends up with the corrrect key $K$ 
(i.e., $x'^{-1}y'^{-1}x'y'=x^{-1}y^{-1}xy$) {\it if and only if} 
$c_b$ commutes with $c_a$. The only visible way to ensure this 
is to have $x' \in A$ and  $y' \in B$.

 Therefore, it appears that if the adversary chooses to solve 
the conjugacy search problem in the group $G$ to recover $x$ and $y$, 
he will then have to face not only the  membership search problem, 
but also the (decision) membership problem, which may very well be 
algorithmically unsolvable. All this seems to be  pushing the adversary 
toward trying to solve a more difficult version of the conjugacy search problem:

\medskip

 {\it Given a group $G$, a subgroup  $A \le G$, and two 
 elements $g, h \in G$,  find $x \in A$ 
such that $h=x^{-1}gx$, given that at least one such $x$ exists.} 
\medskip

 Finally, we note that what we have said in this section does not 
affect some heuristic attacks on the  Anshel-Anshel-Goldfeld protocol 
suggested by several authors \cite{GKTTV1, GKTTV2,  HT}  because these 
attacks, which 
use ``neighbourhood search" type (in a group-theoretic context   
also  called  ``length based") heuristic algorithms, are targeted, 
by design, at finding a solution of a given equation (or a system of 
equations)   as a word in given elements.  
The point that we make in this section is that even if a fast 
(polynomial-time) {\it deterministic} algorithm is found for solving 
the conjugacy search problem in, say, braid groups, this will not 
be sufficient to break the Anshel-Anshel-Goldfeld protocol 
{\it by a deterministic attack}. As for heuristic attacks, their 
limitations are explained in \cite{Shpilrain}.

\end{document}